\documentclass[11pt]{amsart}

\title{Bergman Complexes, Coxeter Arrangements, and Graph Associahedra}

\author{Federico Ardila, Victor Reiner, and Lauren Williams}

\usepackage{pstricks, pst-node, amssymb}
\psset{unit=1pt, arrowsize=4pt, linewidth=.7pt}
\psset{linecolor=blue}
\newgray{grayish}{.90}
\newrgbcolor{embgreen}{0 .5 0}

\def\vblack(#1, #2)#3{\cnode*[linecolor=black](#1, #2){3}{#3}}
\def\vwhite(#1,#2)#3{\cnode[linecolor=black,fillcolor=white,fillstyle=solid](#1,#2){3}{#3}}
\countdef\x=23
\countdef\y=24
\countdef\z=25
\countdef\t=26

\def\tbox(#1,#2)#3{
\x=#1 \y=#2
\multiply\x by 12
\multiply\y by 12
\z=\x \t=\y
\advance\z by 12
\advance\t by 12
\psline(\x,\y)(\x,\t)(\z,\t)(\z,\y)(\x,\y)
\advance\x by 6
\advance\y by 6
\rput(\x,\y){{\bf #3}}}

\newenvironment{caution}{\vspace{2 ex}{\noindent{\bf
Caution.}}}{\vspace{2 ex}}

\usepackage{amsmath, amsthm, amssymb, amsbsy}
\usepackage{amsfonts, latexsym,amscd}
\usepackage[mathscr]{eucal}
\usepackage{epsfig}
\newtheorem{theoremB}{Theorem}
\newtheorem{theorem}{Theorem}[section]

\newtheorem{proposition}[theorem]{Proposition}
\newtheorem{lemma}[theorem]{Lemma}
\newtheorem{def_theorem}[theorem]{Definition/ Theorem}
\newtheorem{definition}[theorem]{Definition}
\newtheorem{example}[theorem]{Example}
\newtheorem{corollary}[theorem]{Corollary}

\newtheorem{remark}[theorem]{Remark}

\newtheorem{observation}[theorem]{Observation}
\usepackage{amsfonts}
\usepackage{amssymb}
\usepackage{amsmath}
\newcommand{\A}{\mathcal{A}}
\newcommand{\R}{\mathbb{R}}
\newcommand{\w}{\omega}
\newcommand{\B}{{\mathcal{B}}}

\newcommand{\FF}{{\mathcal{F}}}

\newcommand{\N}{\mathcal N}

\renewcommand{\L}{{\mathcal L}}
\newcommand{\F}{{\mathcal F}}

\DeclareMathOperator{\depth}{dp}
 \DeclareMathOperator{\cl}{cl}
\DeclareMathOperator{\supp}{supp}

\DeclareMathOperator{\Aut}{Aut}

\DeclareMathOperator{\circuit}{circ}
\DeclareMathOperator{\arbitrary}{arb}
\DeclareMathOperator{\parabolic}{par}
\DeclareMathOperator{\finiteparabolic}{fin par}
\DeclareMathOperator{\interior}{int}

\newcommand{\thmrefer}[1]{\renewcommand\thetheorem
  {\protect\ref{#1}}\addtocounter{theorem}{-1}}

\begin{document}
\maketitle

\section*{Abstract}
Tropical varieties play an important role in algebraic geometry.
The Bergman complex $\B(M)$ and the positive Bergman complex
$\B^+(M)$ of an oriented matroid $M$ generalize to matroids the
notions of the tropical variety and positive tropical variety
associated to a linear ideal. Our main result is that if $\A$ is a
Coxeter arrangement of type $\Phi$ with corresponding
oriented matroid $M_{\Phi}$, then $\B^+(M_{\Phi})$ is dual to the
graph associahedron of type $\Phi$, and $\B(M_{\Phi})$ equals the
nested set complex of $\A$.
In addition, we prove that for any 
orientable matroid $M$, one can find $|\mu(M)|$ different reorientations
of $M$ such that the corresponding positive Bergman complexes cover
$\B(M)$, where $\mu(M)$ denotes the M\"obius function of the lattice of
flats of $M$.  

\section{Introduction}
%

In this paper we study the Bergman complex and the positive
Bergman complex of a Coxeter arrangement.  We relate them to
the nested set complexes that arise in De~Concini and Procesi's
wonderful arrangement models \cite{DP, FK}, and to the graph
associahedra introduced by Carr and Devadoss \cite{CD}, by
Davis, Januszkiewicz,
and Scott \cite{Davis}, and by Postnikov \cite{Postnikov}.

The \emph{Bergman complex} of a matroid is a pure polyhedral
complex which can be associated to any matroid.
It was first defined by Sturmfels \cite{Sturmfels} in
order to generalize to matroids the notion of a tropical variety
associated to a linear ideal.
The Bergman complex  can be described
in terms of the lattice of flats of the matroid, and is homotopy
equivalent to a wedge of spheres, as shown by Ardila and Klivans
\cite{Berg}.

The \emph{positive Bergman complex} $\B^+(M)$
of an oriented matroid $M$ is
a subcomplex of the Bergman complex of the underlying unoriented matroid
$\underline{M}$.
It generalizes to
oriented matroids the notion of the positive tropical variety
associated to a linear ideal.  $\B^+(M)$
depends on a choice of acyclic
orientation of $M$, and as one varies this acyclic orientation,
one gets a covering of the Bergman complex of $\underline{M}$; we
will prove this in Section \ref{2}.
The positive Bergman complex
can be described in terms of the Las Vergnas
face lattice of $M$ and it is homeomorphic to a sphere, as shown
by Ardila, Klivans, and Williams \cite{PosBerg}.

\emph{Graph associahedra} are polytopes which generalize the
associahedron, which were discovered independently by Carr and
Devadoss \cite{CD}, by Davis, Januszkiewicz, and Scott \cite{Davis},
and by Postnikov \cite{Postnikov}.  
There is an intrinsic tiling by associahedra of the
Deligne-Knudsen-Mumford compactification of the real moduli space
of curves $\overline{M_0^n (\R)}$, a space which is related to the
Coxeter complex of type $A$. The motivation for Carr and Devadoss'
work was the desire to generalize this phenomenon to all simplicial
Coxeter
systems.

Let $\A_{\Phi}$ be the Coxeter arrangement corresponding to
the (possibly infinite, possibly non-crystallographic) root system $\Phi$
associated to a Coxeter system $(W,S)$ with diagram $\Gamma$;
see Section~\ref{Cox} below.
Choose a region $R$ of the arrangement, and let $M_{\Phi}$ be the
oriented matroid associated to $\A_{\Phi}$ and $R$. In this paper
we prove:
\begin{theorem}\label{Theorem1}
The positive Bergman complex $\B^+(M_{\Phi})$ of the
arrangement $\A_{\Phi}$ is dual to the graph associahedron
$P(\Gamma)$.
\end{theorem}

In particular, the cellular sphere $\B^+(M_{\Phi})$ is actually
a simplicial sphere, and a {\it flag} (or {\it clique}) complex.

This result is also related to the {\it wonderful model of a hyperplane
arrangement} and to {\it nested set complexes}.
The wonderful model of a hyperplane arrangement is obtained
by blowing up the non-normal crossings of the arrangement, leaving
its complement unchanged. De~Concini and Procesi \cite{DP}
introduced this model in order to study the topology of this
complement. They showed that the \emph{nested sets} of the
arrangement encode the underlying combinatorics.
Feichtner and Kozlov \cite{FK} gave an abstract notion of the
\emph{nested set complex} for any meet-semilattice, and Feichtner
and M\"{u}ller \cite{FM} studied its topology.  Recently,
Feichtner and Sturmfels \cite{FS} studied the relation between the Bergman
complex and the nested set complexes (see Section \ref{Nested} below).

In this paper we also prove for {\it finite} root systems $\Phi$:
\begin{theorem} \label{Theorem2}
The Bergman complex $\B(M_{\Phi})$ of
$\A_{\Phi}$ equals its nested set complex.
\end{theorem}

In particular, the cell complex $\B(M_{\Phi})$ is actually a simplicial complex.

\section{The Bergman complex and the positive
Bergman complex}
 \label{2}

Our goal in this section is to explain the notions of
the Bergman complex of a matroid
and the positive Bergman complex of an oriented matroid
which were studied in \cite{Berg} and \cite{PosBerg}.
In order to do so we must review  a certain operation on matroids
and oriented matroids.

\begin{definition} \rm \ \\
Let $M$ be a matroid or oriented matroid
of rank $r$ on the ground set $[n]$, and let $\w \in \R^n$.  Regard
$\w$ as a weight function on $M$, so that the weight of a basis
$B = \{b_1, \dots , b_r \}$ of $M$ is given by
$\w_B = \w_{b_1} + \w_{b_2} + \dots + \w_{b_r}$.  Let
$B_\w$ be the collection of bases of $M$ having minimum $\w$-weight.
(If $M$ is oriented, then bases in $B_\w$ inherit orientations from
bases of $M$.)
This collection is itself the set of bases of a matroid (or oriented matroid)
which we call $M_\w$.
\end{definition}

It is not obvious that $M_\w$ is  well-defined. However, when $M$
is an unoriented matroid, we can see this by considering the
matroid polytope of $M$: the face that minimizes the linear
functional $\w$ is precisely the matroid polytope of $M_\w$.  For
a proof that $M_\w$ is well-defined when $M$ is oriented, see
\cite{PosBerg}.

Notice that $M_\w$ will not change if we translate $\w$ or scale it
by a positive constant.  We can therefore restrict our
attention to the sphere
$$
S^{n-2} := \{\,\w \in \R^n \,\, : \,\, \w_1 +
\cdots + \w_n = 0 \, , \, \w_1^2 + \cdots + \w_n^2 = 1 \}.
$$
The Bergman complex of $M$ will be a certain subset of this sphere.

The matroid $M_\w$ depends only on a certain flag associated to
$\w$.

\begin{definition} \rm \ \\
Given $\w\in \R^n$, let $\FF(\w)$ denote the unique flag of
subsets
\begin{equation}
\label{flag-of-flats}
\emptyset = F_0 \subset F_1 \subset \dots \subset F_k
\subset F_{k+1}=[n]
\end{equation}
such that $\w$ is constant on each set $F_i \setminus
F_{i-1}$ and satisfies $\w |_{F_i \setminus F_{i-1}} <\w
|_{F_{i+1}\setminus F_{i}}$. We call $\FF(\w)$ the {\it flag} of $\w$, and
we say that the {\it weight class} of $\w$ or of the flag $\FF$ is
the set of vectors $\nu$ such that $\FF(\nu) = \FF$.
\end{definition}

It is shown in \cite{Berg} that $M_\w$ depends only on the flag
$\FF:=\FF(\w)$; specifically
\begin{equation}
\label{matroid-from-flag}
M_\w = \bigoplus_{i=1}^{k+1} F_i / F_{i-1}
\end{equation}
where $F_i/ F_{i-1}$ is obtained from the matroid restriction of
$M$ to $F_i$ by quotienting out the flat $F_{i-1}$.  Hence
we also refer to this oriented matroid $M_\w$ as $M_\FF$.

\begin{def_theorem}\cite{Berg}\label{BergTheorem}
The {\it Bergman complex} of a matroid $M$ on the ground set $[n]$ is
the set
\begin{align*}
{\B}(M)
   & = \{ \w \in S^{n-2} \,\, : \,\, M_{\FF(\w)} \,\,
\text{has no loops} \}\\
   & = \{ \w \in S^{n-2} \,\, : \,\, \FF(\w) \,\,
\text{is a flag of flats of} \,\,M \}
\end{align*}
\end{def_theorem}

Since the matroid $M_{\w}$ depends only on the weight class that
$\w$ is in, the Bergman complex of $M$ is the disjoint union of
the weight classes of flags $\FF$ such that $M_{\FF}$ has no
loops. We say that the weight class of a flag $\FF$ is
\emph{valid} for $M$ if $M_{\FF}$ has no loops.

There are two polyhedral subdivisions of $\B(M)$, one of which is
clearly finer than the other.

\begin{definition} \rm \ \\
\label{coarse-fine}
The \emph{fine subdivision} of $\B(M)$ is the subdivision of
$\B(M)$ into valid weight classes: two vectors $\omega$ and $\nu$ of
$\B(M)$ are in the same class if and only if $\FF(\omega) = \FF(\nu)$.
The \emph{coarse subdivision} of $\B(M)$ is the subdivision of
$\B(M)$ into $M_{\w}$-equivalence classes: two vectors $\omega$ and $\nu$
of $\B(M)$ are in the same class if and only if $M_\omega = M_\nu$.
We call these equivalence classes {\it fine cells} and {\it coarse cells};
however, by default, any reference to a {\it cell} of $\B(M)$ will
refer to a coarse cell.
\end{definition}

The fine subdivision
gives the following corollary of Theorem \ref{BergTheorem}.

\begin{corollary}\cite{Berg}\label{BergTheorem2}
Let $M$ be a matroid of rank $r$.  
The fine subdivision of the Bergman complex
$\B(M)$ is a geometric realization of
$\Delta(L_M-\{\,\hat{0}\,,\hat{1}\,\}\,)$, the order complex of
the proper part of the lattice of flats of $M$. It follows that
$\B(M)$ is homotopy equivalent to a wedge of $|\mu(M)|$ spheres of dimension $r-2$, 
where 
$\mu(M)$ denotes the Mobius function from the bottom to the top element in 
$L_M$.
\end{corollary}

There are positive analogues of all of the above definitions and
theorems.  First we must give the definition of \emph{positive
covectors} and \emph{positive flats}.

\begin{definition} \rm \ \\
Let $M$ be an acyclic oriented matroid on the ground set $[n]$. We
say that a covector $v \in \{+, -, 0\}^n$ of $M$ is
\emph{positive} if each of its entries is $+$ or $0$.  We say that
a flat of $M$ is \emph{positive} if it is the $0$-set of a
positive covector.
\end{definition}

\begin{observation}\label{posflats}
If $M$ is the acyclic oriented matroid corresponding to a
hyperplane arrangement $\A$ whose orientation is determined by
a choice of region $R$, then the
positive flats are in correspondence with the faces of $R$.
In this case we will also say that the flats which are positive are
``positive with respect to $R$.''
\end{observation}

\begin{figure}[h]
\centering
\includegraphics[height=2in]{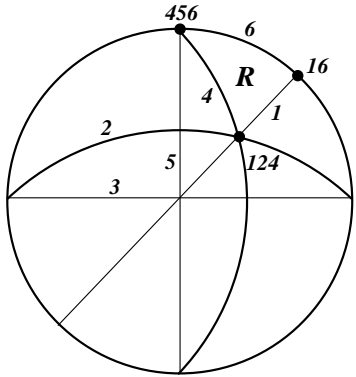}
\caption{The braid arrangement $A_3$.} \label{braid}
\end{figure}

For example, consider the braid arrangement $A_3$, consisting of
the six hyperplanes $x_i=x_j, 1 \leq i<j\leq 4$ in $\R^4$. Figure
\ref{braid} illustrates this arrangement, when intersected with
the hyperplane $x_4=0$ and the sphere $x_1^2+x_2^2+x_3^2=1$. Let
$R$ be the region specified by the inequalities $x_1 \geq x_2 \geq
x_3 \geq x_4$, and let $M_{A_3}$ be the oriented matroid
corresponding to the arrangement $A_3$ and the region $R$. Then
the positive flats are $\emptyset, 1, 4, 6, 124, 16, 456$ and
$123456$.

%
%
%

The positive Bergman complex counterpart to Definition/Theorem
\ref{BergTheorem} is the following.

\begin{def_theorem}\cite{PosBerg}\label{posBerg}
The {\it positive Bergman complex} of $M$ is
\begin{align*}
{\B^+}(M)
   & = \{ \w \in S^{n-2} \,\, : \,\, M_{\FF(\w)} \,\,
\text{is acyclic} \}\\
   & = \{ \w \in S^{n-2} \,\, : \,\, \FF(\w) \,\,
\text{is a flag of positive flats of} \,\,M \}
\end{align*}
\end{def_theorem}

Within each equivalence class of the coarse subdivision of
$\B(M)$, the vectors $\w$ give rise to the same unoriented $M_{\w}$. Since
the orientation of $M_{\w}$ is inherited from that of $M$, they
also give rise to the same oriented matroid $M_{\w}$. Therefore each
coarse cell of $\B(M)$ is either completely contained in or
disjoint from $\B^+(M)$. Thus $\B^+(M)$ inherits the coarse and
the fine subdivisions from $\B(M)$, and each subdivision of
$\B^+(M)$ is a subcomplex of the corresponding subdivision of
$\B(M)$.

Recall that the {\it Las Vergnas face lattice} $\FF_{\ell v}(M)$
is the lattice of positive flats of $M$, ordered by containment.
Note that the lattice of positive flats of the oriented matroid
$M$ sits inside $L_M$, the lattice of flats of $M$. By Observation
\ref{posflats}, if $M$ is the oriented matroid of the arrangement
$\A$ and the region $R$, then $\FF_{\ell v}(M)$ is the face poset
of $R$.

\begin{corollary}\cite{PosBerg} \label{LasVergnas}
Let $M$ be an oriented matroid of rank $r$. Then the fine subdivision of
$\B^+(M)$ is a geometric realization of $\Delta(\FF_{\ell v}(M)
-\{\,\hat{0}\,,\hat{1}\,\}\,)$, the order complex of the proper
part of the Las Vergnas face lattice of $M$.  It follows that
the positive Bergman complex of an oriented matroid is
homeomorphic to an $(r-2)$-sphere.
\end{corollary}


\begin{example}\label{braid-example} \rm
Let $M$ be the oriented matroid from Figure \ref{braid}. The
positive flats of $M$ are $\{\emptyset, 1, 4, 6, 16, 124, 456,
123456\}$. The lattice of positive flats of $M$ is shown in
bold in Figure
\ref{PosFlats}, within the lattice of flats of $M$.


\begin{figure}[h]
\centering
\includegraphics[height=2in ]{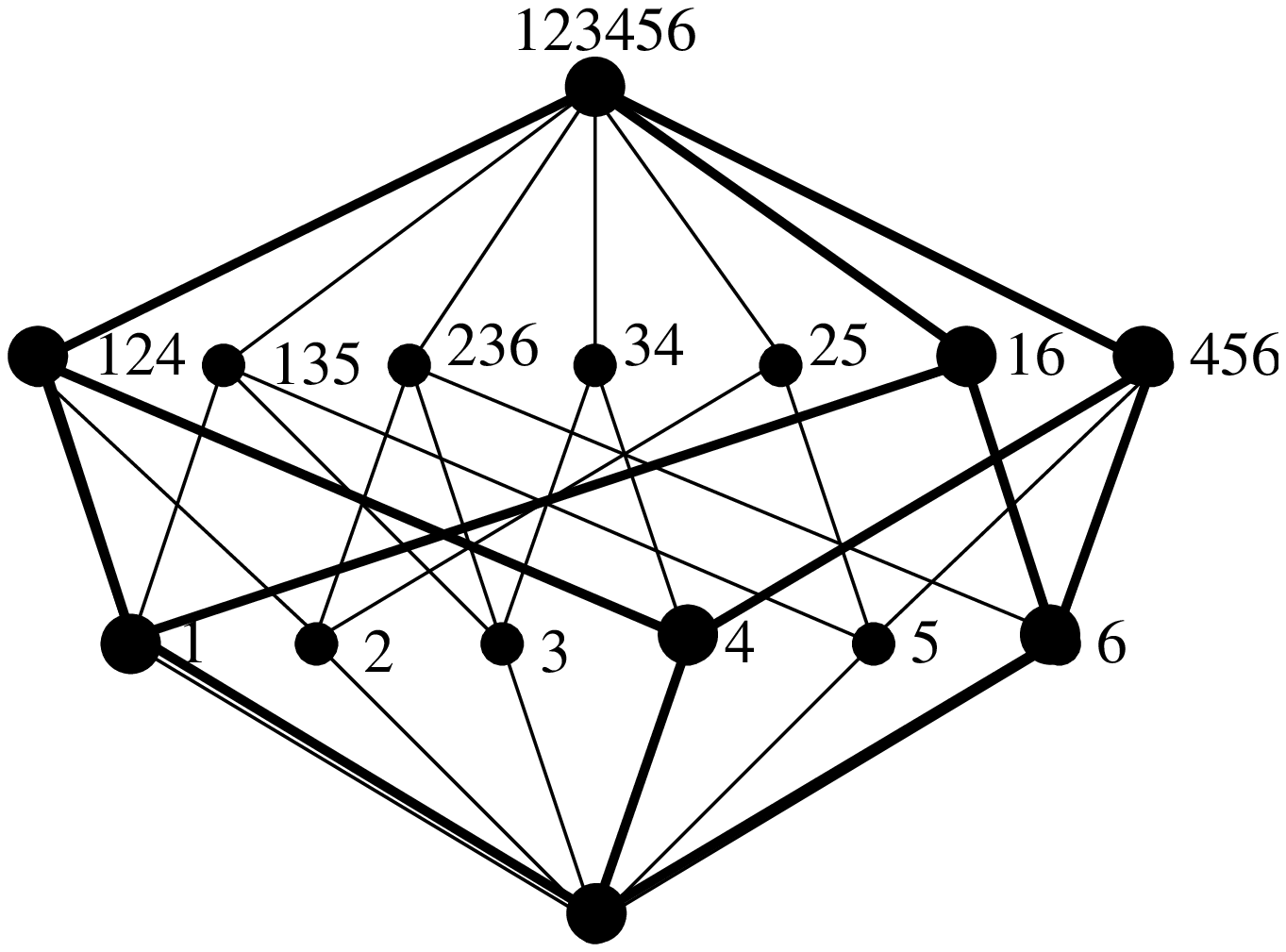}
\caption{The lattice of positive flats within the lattice of flats.}
\label{PosFlats}
\end{figure}
\end{example}

\section{Further theory of Bergman and positive Bergman complexes}
\label{further-Bergman-theory}

  This section develops some further theory of Bergman complexes in
the setting of both unoriented and oriented matroids.  These results will
be used later, in the proofs of Theorems~\ref{Theorem1} and \ref{Theorem2},
but are also of independent interest.

\subsection{Covering the Bergman complex with positive Bergman complexes}

We know that, for any acyclic orientation of a matroid $M$ of rank
$r$, the
corresponding positive Bergman complex is homeomorphic to
an $(r-2)$-sphere, while the Bergman complex $\B(M)$ of the
(unoriented) matroid $M$ is homotopy equivalent to $|\mu(M)|$
such $(r-2)$-spheres. In fact, as we vary the acyclic reorientations (that is, the {\it topes} 
or maximal covectors) of $M$, the corresponding positive Bergman complexes cover 
$\B(M)$.   The first goal of this section is to give a polyhedral realization of this
statement: for any orientable matroid $M$, we exhibit $|\mu(M)|$
reorientations of $M$ whose positive Bergman complexes cover
$\B(M)$.

The motivating example is the matroid $M$ of a real central
hyperplane arrangement $\A$. Let $H$ be an affine hyperplane which
is generic with respect to $\A$. Consider the regions of $\A$
which have a non-empty and bounded intersection with $H$; we will see that there are $|\mu(M)|$ of them. We claim that the positive Bergman complexes
corresponding to these regions cover the Bergman complex of $M$.

For general oriented matroids, one can mimic the previous
construction. Let $M=(E, \L)$ be an oriented matroid of rank $r$
with ground set $E$ and collection of covectors $\L \subseteq
\{+,-,0\}^E$. Let $\L_g \subseteq \{+,-,0\}^{E \cup \{g\}}$ be an
extension of $\L$ by a generic element $g \notin E$; this means
that $g$ is not in the closure of any set $A \subset E$ with $r(A)
< r$. Let $N$ be the affine oriented matroid $N = (E \cup \{g\}, \L_g,
g)$ with distinguished element $g$ (which is not a loop).

Let $\L_g^+ = \{X \in \L_g \, | \, X_g = +\}$ and let
$\widehat{\L}_g^+ = \L_g^+ \cup \{\widehat{0}, \widehat{1}\}$ be
the {\it affine face lattice} of $N$. Let $\L_g^{++} = \{X \in \L_g^+ \,
| \, (\L_g)_{\leq X} \subseteq \widehat{\L}_g^+\}$ be the 
{\it bounded complex} of $N$. The maximal
elements of $\L_g^+$ are the topes, and the maximal elements of
$\L_g^{++}$ are the bounded topes with respect to $g$. Each
bounded tope of $\L_g$ with respect to $g$ determines a tope of
$\L$ by deletion of $g$; we call these the \emph{bounded topes of
$\L$ with respect to $g$}, and write $B^{++}$ for the set of such
topes.

In the realizable case, $M$ is the oriented matroid of $\A$ with
respect to a chosen region. Instead of the generic affine
hyperplane $H$, we consider the translate $g$ of $H$ through the
origin, declaring $H$ to be on the positive side of $g$. The
arrangement $\A \cup \{g\}$ determines the generic extension $N$ of
$M$. The faces in $\L_g^{+}$ correspond to the faces of $\A$ that
intersect $H$, and the faces in $\L_g^{++}$ correspond to the
faces of $\A$ that have a non-empty and bounded intersection with
$H$. The topes in $B^{++}$ are in one-to-one correspondence with
the bounded regions of the arrangement $\A \cup H$.

The \emph{beta 
invariant}  $\beta(N)$ of a matroid $N$ is given by $\beta(N) = 
(-1)^{r(N)} \sum \mu_N(X) r(X)$, summing over all flats $X$ of 
$N$. Here $\mu_N$ denotes the M\"obius function of the lattice 
$L_N$.

\begin{proposition}\cite{GZ, LasVergnas}
An affine oriented matroid $N$ with distinguished element $g$ 
has exactly $\beta(N)$ topes that are bounded with respect to $g$.
\end{proposition}

\begin{lemma}(cf. \cite[Theorem 3.2]{GZ})
If $N$ is a generic extension of $M$ by $g$, then 
$$
\beta(N) = (-1)^{r(N)}\mu(M).$$
\end{lemma}

\begin{proof}
For any nonloop, noncoloop element $g$ in a matroid $N$, one has \cite[Theorem 7.3.2(c)]{Zaslavksy}
\begin{equation}
\label{beta-recurrence}
\beta(N) = \beta(N-g) + \beta(N/g).
\end{equation}
Since $g$ is generic, the lattice $L_{N/g}$ is simply the truncation 
of $L_{N-g}$ in which one removes the entire rank $r-1$ (but keeps 
the element $\widehat{1}$), where $r= r(N)$.  
Hence starting with equation \eqref{beta-recurrence}, 
one has on the right-hand side two sums of the quantities 
$\mu(X) r(X)$ with $X$ ranging over the two lattices $L_{N-g}$ and 
$L_{N/g}$, with opposite signs
in front of the two sums because the ranks of $N-g$ and $N/g$ differ by one.
Thus the terms with $X$ of rank at most $r-2$ all cancel, and one
is left with the terms of rank at least $r-1$ in the two sums:
\begin{eqnarray*}
(-1)^r \beta(N) & = &  r \cdot \mu(N-g)
    + (r-1)\left( \sum_{X \mathrm{\,\, of \,\, rank\,\, } 
      r-1 \mathrm{\,\, in\,\, }L_{N-g}} \mu(X)\right)
          - (r-1) \mu(N/g) \\
 & = & r \cdot \mu(N-g)
    + (r-1)\left( \mu(N/g) - \mu(N-g)\right)
          - (r-1) \mu(N/g) \\
 & = &\mu(N-g),  
\end{eqnarray*}
as we wished to show. \end{proof}

\begin{theorem}\label{posBergmanscoverBergman}
Let $M$ be an oriented matroid, and $N$ an extension by a generic
element $g$. Let $T_1, \ldots, T_{|\mu(M)|}$ be the bounded
topes in $M$ with respect to $g$. Then the $|\mu(M)|$ positive
Bergman complexes corresponding to the $T_i$s cover the Bergman
complex of the unoriented matroid $\underline{M}$.
\end{theorem}

\begin{proof}
There is no harm in assuming that $M$
is simple and loop and coloop-free. In view of Definition/Theorem
\ref{posBerg}, it
suffices to show that, for any flag of flats $\F = \{\emptyset
\subset F_1 \subset \cdots \subset F_{r-1} \subset E\}$, we can
find a tope $T_i$ with $1 \leq i \leq |\mu(M)|$ such that all
the $F_j$s are positive with respect to $T_i$. This means that
$T_i$ has a flag of subfaces (covectors with some entries of $T_i$ replaced by zeroes) $X_1 > \ldots >
X_{r-1}$ such that $X_i$ spans $F_i$.

We proceed by induction, where the base case is trivial. Now
consider the rank $r-1$ oriented matroid $M/F_1$ (which, in the
realizable case, corresponds to the arrangement that $\A$
determines on the hyperplane $F_1$). The set $F_1$ is also a flat
in $N$, and $N/F_1$ is also a generic extension of $M/F_1$ by $g$.
Consider the flag of flats 
$$
\F' = \{\emptyset \subset F_2-F_1
\subset \cdots \subset F_{r-1}-F_1 \subset E-F_1\}
$$ 
of $M/F_1$. By
the induction hypothesis, we can find a tope $T'$, bounded in
$M/F_1$ with respect to $g$, which has a flag of faces $T'=Y_1
> Y_2 > \ldots > Y_{r-1}$ such that $Y_i$ spans $F_i-F_1$ in $M/F_1$.

Since $M$ is simple, the flat $F_1$ consists of a single element; call it $e$.
Then the covector $Y_i$ of $M/F_1$ comes from the covector $X_i$ of $M$
which is identical to $Y_i$, except for the extra entry
$(X_i)_e=0$. Clearly $X_1 > \cdots > X_{r-1}$ and $X_i$ spans
$F_i$.

Since $e$ is not a loop of $M$, some covectors $Z^+$ and $Z^-$ of
$M$ have $Z^+_e=+$ and $Z^-_e=-$. The topes $T^+ = X_1 \circ Z^+$
and $T^- = X_1 \circ Z^-$ are identical to $Y_1=T'$, except for
the extra entries $(T^+)_e = +$ and $(T^-)_e = -$. The $X_i$s are
faces of both $T^+$ and $T^-$, so it remains to show that at least
one of $T^+$ and $T^-$ is bounded with respect to $g$ in $M$.

Suppose this is not the case.
Then we can find non-zero covectors $A \leq T^+\cup \{g\}$ and $B \leq
T^- \cup \{g\}$ of $N$ such that $A_g=B_g=0$. If $A_e =0$, then $A-e$
would be a non-zero covector in $N/F_1$, smaller than $T'\cup \{g\}$
and satisfying $(A-e)_g=0$; this would contradict the boundedness
of $T'$ in $M/F_1$ with respect to $g$. Therefore we have $A_e =
+$ and, similarly, $B_e = -$.

Consider now the covectors $A$ and $B$ of $N$ and their separator
$e$. In fact, $e$ is the only separator of $A$ and $B$, because 
$$
\begin{aligned}
A \leq T^+ \cup \{g\} &= T' \cup \{e\} \cup \{g\}\\
B \leq T^- \cup \{g\} &= T' \cup \{\overline{e}\} \cup \{g\}.\\
\end{aligned}
$$ 
Using the covector axiom (L3)
\cite[Theorem 4.1.1]{RedBook}, we will find a covector $C$ of $N$
such that $C_e=0$ and $C_f = (A \circ B)_f = (B \circ A)_f$ for
all $f \neq e$. Thus $C-e$ is a covector of $N/F_1$ which is
smaller than $T' \cup \{g\}$ and satisfies $(C-e)_g=0$. This
contradicts the boundedness of $T'$ in $M/F_1$ with respect to
$g$, unless $C-e=0$. But if this were the case, then $e$ would be a
coloop of $N$. In the presence of the generic element $g$, this is
impossible: $E-e$ has corank at most $1$ in $M$, so $( E- e )\cup \{g\}$
is spanning in $N$, without containing $e$. This completes the
proof.
\end{proof}

Theorem \ref{posBergmanscoverBergman} is closely related to recent
work of Bj\"orner and Wachs. In \cite{BW}, they construct a basis
for the homology of the geometric lattice of an orientable matroid
$M$, which is indexed by the bounded topes of $M$ with respect to
an extension by a generic element $g$.

\subsection{The forest of a flag, and coarse cells in the Bergman complex}

Recall from Definition~\ref{coarse-fine} that the Bergman complex
$\B(M)$ has two subdivisions into cells.
Its fine subdivision has cells indexed by all flags $\FF$ of flats of $M$.  These
fine cells then group themselves into the cells of the coarse subdivision, according to their 
associated matroids $M_\FF$.
It turns out that one can always determine $M_\FF$, and 
hence the coarse cell to which a flag $\FF$ corresponds, based on 
a certain labelled forest $T_\FF$ associated to $\FF$.  These forests also turn
out (see Section~\ref{Nested}) to be closely related to the complex of nested sets\footnote{The material in
this subsection is closely related to results of Feichtner and Sturmfels
\cite[Section 4 and end of Section 3]{FS};  in particular, see our Remark~\ref{FS-remark}.
However, the crucial notion of a {\it circuitous base} 
(Definition~\ref{circuitous-defn}, Proposition~\ref{circuitous-prop}) does not appear in their work, and
we have chosen to explain this subsection in our language so as to keep the paper more self-contained.}.

Recall that the {\it connected components}
of a matroid $M$ are the equivalence classes for the following
equivalence relation on the ground set $E$ of $M$:  say $e \sim e'$ for two elements
$e, e'$ in $E$ whenever they lie in a common circuit of $M$, and then take the transitive
closure of $\sim$.  Recall also
that every connected component is a flat of $M$, and $M$ decomposes (uniquely)
as the direct sum of its connected components.

\begin{definition} \rm \ \\
To each flag $\FF$ of flats of a matroid $M$ indexed as in \eqref{flag-of-flats},
associate a forest $T_\FF$ of rooted trees, in which each vertex $v$ is labelled
by a flat $F(v)$, as follows:
\begin{enumerate}
\item[$\bullet$] For each connected component $F$ of the matroid $M$,
create a rooted
tree (as specified below) and label its root vertex with $F$.
\item[$\bullet$] For each vertex $v$ already created, and already labelled by some flat $F(v)$
which is a connected component of some flat $F_j$ in the flag $\FF$, create
children of $v$ labelled by each of the connected components of $F_{j-1}$ which are contained
{\it properly} in $F(v)$.
\end{enumerate}
\end{definition}

Alternatively, one can construct the forest $T_\FF$ by listing all
the connected components of all the flats in $\FF$, and partially
ordering them by inclusion.

\begin{proposition}
\label{forest-determines-matroid}
For any flag $\FF$ of flats in a matroid $M$, the labelled forest $T_\FF$
determines the matroid $M_\FF$.
\end{proposition}
\begin{proof}
Recall the expression \eqref{matroid-from-flag} for $M_\FF$.
By construction of $T_\FF$, every component of $F_i$ is $F(v)$ for some unique vertex
$v$, and every component of $F_{i-1}$ lying in $F(v)$ is $F(v')$ for some child $v'$ of $v$.
Since quotients commute with direct sums, this gives
\begin{equation}
\label{matroid-from-tree} M_\FF=\bigoplus_{\text{ vertices }v
\text{ of }T_\FF}
      \left( F(v) \slash \bigoplus_{\text{ children }v'\text{ of }v} F(v') \right) .
\end{equation}
\end{proof}

In general, the converse of this proposition does {\it not} hold;
one can have $M_{\FF}=M_{\FF'}$ without $T_\FF=T_{\FF'}$.  For
example (cf. \cite[Example 1.2]{FS}), in the matroid $M$ on ground
set $E=\{1,2,3,4,5\}$ having rank $3$ and circuits $\{123, 145,
2345\}$, the two flags
$$
\begin{aligned}
\FF &:=( \emptyset \subset 1 \subset 123 \subset 12345 ) \\
\FF'& :=( \emptyset \subset 1 \subset 145 \subset 12345 )
\end{aligned}
$$
exhibit this possibility.

However, there is at least one nice hypothesis that allows one to
reconstruct $T_\FF$ from $M_\FF$.  Given a base $B$ of a matroid $M$ on ground
set $E$, and any element $e \in E \backslash B$, there is a unique circuit
of $M$ contained in $B \cup \{e\}$, called the {\it basic circuit} $\circuit(B,e)$.
Note that the flat spanned by $\circuit(B,e)$ will always be a connected flat.

\begin{definition} \rm \ \\
\label{circuitous-defn}
Say that a base $B$ of a matroid $M$ is {\it circuitous} if every
connected flat spanned by a subset of $B$ is spanned by the basic
circuit $\circuit(B,e)$ for some $e \in E \backslash B$.
\end{definition}

Note that the basic circuit $\circuit(B,e)$ spanning the connected
flat $F$ must be $(F \cap B) \cup \{e\}$. Before we state our
proposition, we prove two useful lemmas.

\begin{lemma}
\label{lemma1} Let $F$ be a flat in a matroid $M$, spanned by some
independent set $I$.  Then every connected component of F is
spanned by some subset of $I$, namely, by the intersection of that
component with I.
\end{lemma}

\begin{proof}
Let $r$ denote the rank function for $M$, and let
$F$ have components $F_1,\ldots,F_t$. Then
\[
\sum_{i} r(F_i) = r(F)= |I| = \sum_{i} | F_i \cap I | = \sum_{i}
r(F_i \cap I) \leq \sum_{i} r(F_i),
\]
which means we must have an equality for each $i$: $r(F_i \cap I)
= r(F_i)$. In other words, $F_i \cap I$ spans $F_i$.
\end{proof}

Given a subset $A \subset E$ of the ground set of a matroid,
let $\cl(A)$ denote its {\it closure}, that is, the flat spanned
by $A$.

\begin{lemma}
\label{lemma2} Let $F \subset G$ be flats of a matroid that are
spanned by subsets of a circuitous base $B$. If $G$ is connected,
then $G/F$ is also connected.
\end{lemma}

\begin{proof}
Let $I_F = F \cap B$ and $I_G = G \cap B$; these are bases for $F$
and $G$, respectively. Also, $I_F \subset I_G$, and $I_G - I_F$ is
a base for the quotient $G/F$. Since $G$ is a connected flat
spanned by a subset of the circuitous base $B$, there exists $e$
in $G - B$ such that $\cl(\circuit(B,e))=G$, and $\circuit(B,e) =
I_G \cup \{e\}$.

We now claim that 
$$
\circuit_{G/F}( I_G - I_F, e ) = I_G - I_F \cup \{e\}.
$$ 
We need to check that $I_G - I_F \cup \{e\} - \{g\}$ is independent in
$G/F$ for any $g \in I_G - I_F \cup \{e\}$. Since $I_F$ is a basis of
$F$, this follows from the fact that $I_G \cup \{e\} - \{g\}$ is
independent in $G$. We conclude by observing that $G/F$ is the
flat spanned by $\circuit(I_G - I_F, e)$, so it is connected.
\end{proof}

\begin{proposition}
\label{circuitous-prop}
Let $B$ be a circuitous base of a matroid $M$.  Then for any two flags $\FF, \FF'$
of flats spanned by subsets of $B$,
one has $M_\FF = M_{\FF'}$ if and only if $T_\FF = T_{\FF'}$.
\end{proposition}
\begin{proof}
We start by making two observations about the matroid $M_{\FF}$
and the tree $T_{\FF}$.

First we observe that, under these hypothesis, the expression
\eqref{matroid-from-tree} is actually the decomposition of
$M_{\FF}$ into its irreducible components. By Lemma \ref{lemma1},
the $F(v)$s are connected flats spanned by subsets of $B$. The
direct sums $\oplus_{v '} F(v')$ are also spanned by subsets of
$B$. Lemma \ref{lemma2} then guarantees that $F(v)/\oplus_{v '}
F(v')$ is connected for each vertex $v$ of the tree.

Secondly we show that, among the sets $\cl(\circuit(B,e))$ with
$e$ in $F(v) \backslash \cup F(v')$ and not in $B$, there is a
maximum one under containment, which is precisely $F(v)$.

Take any $e$ in $F(v) \backslash \cup F(v')$ and not in $B$. The
flat $F(v)$ is spanned by a subset $I$ of $B$, and $I \cup \{e\}$ is
dependent. Therefore $\circuit(B,e) \subseteq I \cup \{e\} \subseteq
F(v)$, which implies $\cl(\circuit(B,e)) \subseteq F(v)$.

Now, since $F(v)$ is a connected flat spanned by a subset of $B$,
$F(v) = \cl(\circuit(B,e))$ for some $e \in E \backslash B$.
Clearly $e \in F(v)$. If $e$ was in $F(v')$ for some child $v'$ of
$v$, the argument of the previous paragraph would imply that
$\cl(\circuit(B,e)) \subseteq F(v')$. Therefore $e \in F(v)
\backslash \cup F(v')$.

The two previous observations give us a procedure to recover the
tree $T_\FF$ from the matroid $M_{\FF}$. The first step is to
decompose $M_\FF$ into its connected components $M_1, \ldots,
M_t$, having accompanying ground set decomposition $E = E_1 \sqcup
\cdots \sqcup E_t$. The second step is to recover the flat
corresponding to each $M_i$, as the maximum $\cl(\circuit(B,e))$
with $e \in E_i \backslash B$. The labelled forest $T_\FF$ is
simply the poset of inclusions among these flats.
%
%
%
\end{proof}

It will turn out that the simple roots $\Delta$ of a root system $\Phi$
always form a circuitous base for the associated matroid $M_\Phi$; see
Proposition~\ref{Coxeter-flats}(iii) below.

\begin{remark} \rm
\label{FS-remark}
When the matroid $M$ is connected, the forest $T_\FF$ constructed above is a rooted tree.
It coincides with the tree constructed by Feichtner and Sturmfels in \cite[Proposition 3.1]{FS}
when they choose the minimal building set for their lattice.  In this way,
Proposition~\ref{forest-determines-matroid}
follows from \cite[Theorem 4.4]{FS}.
\end{remark}

\section{Graph associahedra}

Graph associahedra are polytopes which generalize the
associahedron, which were discovered independently by Carr and
Devadoss \cite{CD}, Davis, Januszkiewicz, and Scott \cite{Davis},
and Postnikov \cite{Postnikov}. There is an intrinsic tiling by
associahedra of the Deligne-Knudsen-Mumford compactification of
the real moduli space of curves $\overline{M_0^n (\R)}$, a space
which is related to the Coxeter complex of type $A$. The
motivation for Carr and Devadoss' work was the desire to
generalize this phenomenon to other Coxeter systems.



In order to define graph associahedra, we must introduce the
notions of tubes and tubings.  We follow the presentation of
\cite{CD}.

\begin{definition} \rm \ \\
Let $\Gamma$ be a graph.  A {\it tube} is a proper nonempty set of
nodes of $\Gamma$ whose induced graph is a proper, connected
subgraph of $\Gamma$.  There are three ways that two tubes can
interact on the graph:
\begin{itemize}
\item Tubes are {\it nested} if $t_1 \subset t_2$.
\item Tubes {\it intersect} if $t_1 \cap t_2 \neq \emptyset$ and
  $t_1 \not\subset t_2$ and $t_2 \not\subset t_1$.
\item Tubes are {\it adjacent} if $t_1 \cap t_2 = \emptyset$
  and $t_1 \cup t_2$ is a tube in $\Gamma$.
\end{itemize}

Tubes are {\it compatible} if they do not intersect and they are
not adjacent.  A {\it tubing} $T$ of $\Gamma$ is a set of tubes
of $\Gamma$ such that every pair of tubes in $T$ is compatible.
A $k$-{\it tubing} is a tubing with $k$ tubes.
\end{definition}

Graph-associahedra are defined via a construction which we will now
describe.

\begin{definition} \rm \ \\
Let $\Gamma$ be a graph on $n$ nodes.  Let $\Delta_{\Gamma}$ be the
$n-1$ simplex in which each facet corresponds to a particular node.
Note that each proper subset of nodes of $\Gamma$ corresponds to
a unique face of $\Delta_{\Gamma}$, defined by the intersection of
the faces associated to those nodes.  The empty set corresponds
to the face which is the entire polytope $\Delta_{\Gamma}$.
For a given graph $\Gamma$, truncate faces of $\Delta_{\Gamma}$ which
correspond to $1$-tubings in increasing order of dimension
(i.e. first truncate vertices, then edges, then $2$-faces, \dots).
The resulting polytope $P(\Gamma)$ is the {\it graph associahedron}
of Carr and Devadoss.
\end{definition}

\begin{figure}
\centering
\includegraphics[height=1.5in]{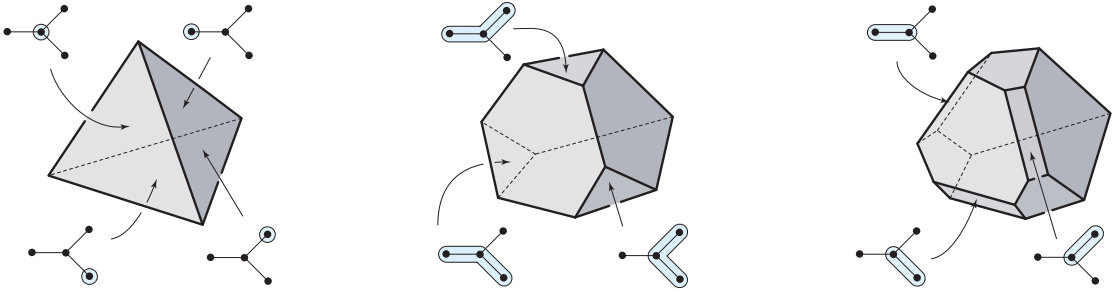}
\caption{$P(D_4)$ \quad \copyright {\tiny Satyan Devadoss}}
\label{D4}
\end{figure}

Figure \ref{D4} illustrates the construction of the
graph associahedron of a Coxeter diagram of type $D_4$. We start
with a simplex, whose four facets correspond to the nodes of
the diagram. In the first step, we truncate three of the vertices,
to obtain the second polytope shown. We then truncate three of the
edges, to obtain the third polytope shown. In the final step, we
truncate the four facets which all correspond to tubes. This step
is not shown in Figure \ref{D4}, since it does not affect the
combinatorial type of the polytope.

\begin{figure}
\centering
\includegraphics[height=1.5in]{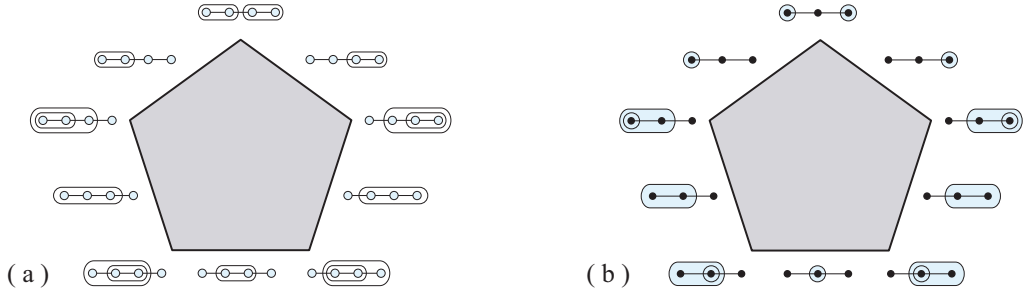}
\caption{The associahedron $A_2$ is the graph associahedron of a $3$-element chain. \quad \copyright {\tiny Satyan Devadoss}}
\label{A3}
\end{figure}

When the graph $\Gamma$ is the $n$-element chain, the
polytope $P(\Gamma)$ is the associahedron $A_{n-1}$. One can see
this by considering an easy bijection between valid tubings and
parenthesizations of a word of length $n-1$, as illustrated in
Figure \ref{A3}.

%

Carr and Devadoss proved that the face poset of $P(\Gamma)$ can
be described in terms of valid tubings.

\begin{theorem}\cite{CD} \label{Face}
The face poset of $P(\Gamma)$ is isomorphic to the set of valid
tubings of $\Gamma$, ordered by reverse containment:
$T < T^{\prime}$ if
$T$ is obtained from $T^{\prime}$ by adding tubes.
\end{theorem}

\begin{corollary}\cite{CD}
When $\Gamma$ is a path with $n-1$ nodes, $P(\Gamma)$ is the
associahedron $A_{n}$ of dimension $n$.  When $\Gamma$ is a cycle with $n-1$
nodes, $P(\Gamma)$ is the cyclohedron $W_n$.
\end{corollary}

\section{Coxeter systems, the Tits cone, and parabolic flats}
\label{Tits-cone}

In this section we review the notion of a Coxeter system $(W,S)$, and
explain two ways of thinking about the associated matroid.  The first
way is to consider the vector configuration of positive roots $\Phi^+$
of the corresponding root system $\Phi$ in $V:= \R^{|S|}$.  The second way
is to
consider a certain arrangement $\A_\Phi$ of hyperplanes in $V^*$
intersecting a
$W$-invariant convex cone known as the Tits cone.
We attempt to give a careful
discussion of the issues that arise when $W$ is infinite,
for example, how to define the Bergman complex,
and what kinds of flats
of the associated matroid are relevant for
Bergman complexes and wonderful compactifications.

A {\it Coxeter system} is a pair $(W,S)$ consisting of a group $W$
and a set of generators $S \subset W$, subject only to relations
of the form $$(s s^{\prime})^{m(s,s^{\prime})} = 1,$$
where $m(s,s)=1$, $m(s,s^{\prime})=m(s^{\prime},s) \geq 2$ for
$s \neq s^{\prime}$ in $S$.  In case no relation occurs for a pair
$(s,s^{\prime})$, we make the convention that $m(s,s^{\prime})=\infty$.
We will always assume that $S$ is finite.

Note that to specify a Coxeter system $(W,S)$, it is enough to draw
the corresponding {\it Coxeter diagram} $\Gamma$: this is a graph
on vertices indexed by elements of $S$, with vertices $s$ and $s^{\prime}$
joined by an edge labelled $m(s,s^{\prime})$ whenever this number
($\infty$ allowed) is at least $3$.

\begin{remark} \rm
In what follows, the reader should note that the {\it positive} Bergman
complex and the graph associahedron associated
with $\Gamma$ will turn out {\it not} to depend on the
edge labels $m(s,s^\prime)$ of $\Gamma$, and only
depend upon the undirected graph underlying $\Gamma$.  However, the Bergman complex
{\it will} turn out to depend upon the edge labels  $m(s,s^\prime)$.
\end{remark}

Although an arbitrary Coxeter system $(W,S)$
need not have a faithful representation
of $W$ as a group generated by orthogonal reflections (for
a positive definite inner product), there exists
a reasonable substitute, called its {\it geometric representation}
\cite[Sec. V.4]{Bourbaki}, \cite[Sec. 5.3, 5.13]{Humphreys}, which we recall here.
Let $V:=\R^{|S|}$ with a basis of {\it simple roots} $\Delta:=\{\alpha_s: s \in S\}$.
Define an $\R$-valued bilinear form $( \cdot,\cdot )$ on $V$ by
$$
( \alpha_s, \alpha_{s'} ):= -\cos \left(\frac{\pi}{m(s,s^\prime)} \right)
$$
and let $s$ act on $V$ by the ``reflection'' that fixes $\alpha_s^\perp$
and negates $\alpha_s$:
$$
s(v):=v - 2(v,\alpha_s) \alpha_s.
$$
This turns out to extend to a faithful representation of $W$ on $V$,
and one defines the {\it root system} $\Phi$ and {\it positive roots} $\Phi^+$ by
$$
\begin{aligned}
\Phi &:=\{w(\alpha_s): w \in W, s \in S\} \\
\Phi^+ & := \{ \alpha \in \Phi: \alpha = \sum_{s \in S} c_s \alpha_s \text{ with }c_s \geq 0 \}
\end{aligned}
$$
It turns out that $\Phi = \Phi^+ \sqcup \Phi^-$ where
$\Phi^-:=-\Phi^+$, and that $W$ will be infinite if and only if $\Phi$ is infinite.

\begin{definition} \rm \ \\
Given a root $\beta \in \Phi$, expressed uniquely in terms of
the simple roots $\Delta$ as  $\beta = \sum_{s \in S} c_s \alpha_s$,
define the {\it support} of $\beta$ (written $\supp \beta$)
to be the vertex-induced subgraph of the Coxeter diagram $\Gamma$
on the set of  vertices $s \in S$ for which  $c_s \neq 0$.
\end{definition}

We will need the following lemma about supports of roots.  It is well-known
when $W$ is finite and crystallographic \cite[No. VI.1.6, Cor. 3]{Bourbaki}, and
a proof of its first assertion for the Coxeter systems associated to Kac-Moody Lie algebras
can be found in \cite[Lemma 1.6]{Kac}; we will need the assertion
in general.

\begin{lemma}\label{KacLemma}
Let $(W,S)$ be an arbitrary Coxeter system with Coxeter graph $\Gamma$.
Then for any root $\beta \in \Phi$ the graph $\supp \beta$ is
connected, and conversely, every connected subgraph $\Gamma'$ of $\Gamma$ occurs
as $\supp \beta$ for some positive root $\beta$.
\end{lemma}

\begin{proof}
For the first assertion, let $\beta$ be a root, which we may assume
is positive without loss of generality.  It is known \cite[Sec. 4.6]{BjornerBrenti}
that there exists a chain of (distinct) positive roots
$$
\alpha=\beta_1 \lessdot \beta_2 \lessdot \cdots \lessdot \beta_k = \beta
$$
in which $\alpha$ is a simple root, and where each relation $\lessdot$ is a covering relation
in what Bj\"orner and Brenti call the {\it root poset}.  
This is the poset on positive roots defined as follows:
$\beta \leq \gamma$ if there exists $s_1, s_2, \dots , s_k \in S$ such that
\begin{enumerate}
\item $\gamma = s_k s_{k-1} \dots s_1 \beta$, and 
\item $\depth(s_i s_{i-1} \dots s_1 \beta) = \depth(\beta) +i$, for all
$1 \leq i \leq k$.
\end{enumerate}
Here, the {\it depth} $\depth$ of a positive root is 
defined to be 
\begin{equation*}
\depth(\beta) = min\{k: w(\beta) \in \Phi^{-} \text{ for some }
w \in W \text{ with } \ell(w) = k \}.
\end{equation*}
In particular, when 
two positive roots
$$
\begin{aligned}
\gamma&=\sum_{t \in S} c_t \alpha_t\\
\gamma'&=\sum_{t \in S} c'_t \alpha_t
\end{aligned}
$$ 
satisfy $\gamma \lessdot \gamma'$, then 
\begin{equation}
\label{gamma'-expression}
\gamma'= s(\gamma) = \gamma - 2(\gamma,\alpha_s) \alpha_s
\end{equation}
for some $s \in S$, and the nonnegative coefficients 
$c_t, c_t'$ satisfy $c_t \leq c_t'$ for all $t \in S$,
so that $\supp \gamma \subseteq \supp \gamma'$; 
see \cite[Corollary 4.6.5]{BjornerBrenti}.
By induction on $k$, it suffices to show that if $\supp(\gamma)$ is connected and
$\gamma \lessdot \gamma'$, then $\supp(\gamma')$ is connected.
From expression \eqref{gamma'-expression} we conclude that 
$\supp\gamma  \subseteq \supp\gamma' \subseteq \supp\gamma \cup \{s\}$.  
Hence either 
\begin{enumerate}
\item[$\bullet$] $\supp\gamma' = \supp\gamma$, which is connected, so we are done, or
\item[$\bullet$] $s \not\in \supp\gamma$ and $\supp\gamma' = \supp\gamma \sqcup \{s\}$.
If $\supp\gamma \sqcup \{s\}$ is connected, we are done.  If not, then $(\gamma,\alpha_s)=0$,
so the expression \eqref{gamma'-expression} forces the contradiction $\gamma'= \gamma$.
\end{enumerate}

For the second assertion, let $\Gamma'$ be a connected subgraph of $\Gamma$,
and we will exhibit a positive root $\gamma'$ with $\supp \gamma' =\Gamma'$ using induction
on the number of vertices of $\Gamma'$.  Let
$s \in S$ be a vertex lying in $\Gamma'$ whose removal leaves a connected subgraph
$\Gamma''=\Gamma - \{s\}$.  By induction there exists a positive root $\gamma$ having
$\supp \gamma  = \Gamma''$, and we claim that $\gamma':=s(\gamma)$ has $\supp(\gamma')=\Gamma'$.
To see this, note that $\gamma=\sum_{t \in \Gamma''} c_t \alpha_t$ with each $c_t > 0$.
Hence
$$
(\alpha_s,\gamma) = \sum_{t \in \Gamma''} c_t ( \alpha_s, \alpha_t ) < 0
$$
since each $( \alpha_s, \alpha_t )$ is nonpositive, and at least one is negative
due to $\Gamma'' \cup \{s\} = \Gamma'$ being connected.  Therefore
the expression \eqref{gamma'-expression} for $\gamma'$ shows that $\supp(\gamma')=\Gamma'$
\end{proof}

We use $M_\Phi$ to denote the matroid
represented by the vector configuration of positive roots
$\Phi^+$ in $V$.  Thus $M_\Phi$ is a matroid of finite rank $r=|S|$,
but has ground set $E=\Phi^+$ of possibly (countably) infinite
cardinality.

\begin{remark} \rm \label{infinite}
When the ground set $E$ is infinite, we need to be careful about how
we define the objects that we are studying: it is no longer
clear what is meant by a weight vector $\w$ or the
bases of minimum $\w$-weight.  Therefore we will not refer to
$M_\w$ in this case; only to
the matroid $M_\FF$ associated to a flag of flats $\FF$.
We will not think of the positive Bergman complex as a subset
of weight vectors (as in Definition/Theorem \ref{posBerg}), but as
a coarsening of the order complex of the lattice of positive
flats (which we can do by Corollary \ref{LasVergnas}).
Although we can similarly consider the Bergman complex  as a coarsening
of the order complex of the lattice of flats, for technical reasons
we will not
deal with the Bergman complex of a matroid with an infinite ground set
in this paper.
\end{remark}

For an arbitrary Coxeter system $(W,S)$, when
one wants to think of the oriented matroid $M_\Phi$ as the
oriented matroid of a hyperplane arrangement $\A_\Phi$ (as opposed to the oriented
matroid of the configuration of vectors $\Phi^+$), one must
work with the contragredient representation $V^*$.  Then $M_\Phi$ is
simply the matroid of the reflecting hyperplanes in $V^*$  for the positive
roots $\Phi^+$.

We now review the Tits cone.  See \cite[Sec. V.4]{Bourbaki},
\cite[Sections 1.15, 5.13]{Humphreys},
\cite[Chapter I]{Brown}, and particularly \cite{Walter} for a very detailed discussion.
Let $\{\delta_s: s \in S\}$ denote the basis for $V^*$ dual to the basis of
simple roots $\Delta$ for $V$.
Then the {\it (closed) fundamental chamber}
$R$ is the nonnegative cone spanned by $\{\delta_s: s \in S\}$ inside $V^*$.
The {\it Tits cone} is the union $T:=\bigcup_{w \in W} w(R)$, a (possibly proper, not
necessarily closed nor polyhedral) convex cone inside $V^*$.  Every positive root
$\alpha \in \Phi^+$ has associated a hyperplane and two half-spaces, $H_\alpha,
H^+_\alpha, H^-_\alpha$ in $V^*$, consisting of those functionals $f \in V^*$
for which $f(\alpha)$ is zero, positive, or negative, respectively.
These hyperplanes and half-spaces decompose the Tits cone\footnote{When
$W$ is infinite, note that only part of the hyperplane or its half-spaces lies inside the Tits cone
$T$.} into cells $\sigma$ that turn out
to be simplicial cones $\sigma$, each of them relatively open within the linear subspace that they span.
The top-dimensional (open) cones
are exactly the images $w(\interior(R))$ as $w$ runs through $W$, where $\interior(R)$ denotes the
{\it interior} of the fundamental chamber $R$.  The tope
(maximal covector) in the  oriented matroid $M_\Phi$ associated to
$w(\interior(R))$ will have the sign $+$ on the roots $\Phi^+ \cap w^{-1}(\Phi^+)$ and
the sign $-$ on the roots $\Phi^+ \cap w^{-1}(\Phi^-)$.
More generally one has the following proposition (see \cite[Section 3.2]{Walter}) relating
an arbitrary cone $\sigma$ to a coset $wW_J$ of a {\it standard parabolic subgroup} $W_J$ (= the subgroup of
$W$ generated by $J$) for some subset $J \subseteq S$

\begin{proposition}
The cones $\sigma$ in the decomposition of $T$ are naturally
in bijection with the cosets $wW_J$
of standard parabolic subgroups, with $wW_J$ determined by the following equality:
$$
wW_J =\{u \in W: u^{-1}(\sigma) \subseteq R\}.
$$

The cone $\sigma$ will then have the following description as an intersection:
for any $u$ in
the coset $wW_J$, one has
$$
\sigma = \bigcap_{s \in J} u(H_{\alpha_s}) \cap \bigcap_{s \in S \backslash J}
u(H^+_{\alpha_s}).
$$
\end{proposition}

As a consequence of this proposition (see \cite[Prop. 3.4]{Walter}),
the linear span of the cone $\sigma$ in $V^*$ is the
hyperplane intersection $\bigcap_{s \in J} H_{w(\alpha_s)}$,
which is the subspace $\left(V^*\right)^{wW_Jw^{-1}}$
fixed by the parabolic subgroup $wW_Jw^{-1}$.

As pointed out in Remark \ref{infinite}, when $W$ (equivalently $\Phi$, or $E=\Phi^+$)
is infinite,
we want to consider the positive
Bergman complex to be a coarsening of the order complex of the lattice
of positive flats.  However, in this situation we have a choice to make,
because there are
three different kinds of flats $F$ of the (oriented) matroid $M_\Phi$,
not all of which are relevant to the De~Concini-Procesi wonderful
compactifications.  These three kinds of flats are
distinguished by how the associated intersection subspace
$$
X_F:=\bigcap_{\alpha \in F} H_\alpha
$$
intersects the Tits cone $T$:

\begin{enumerate}
\item An {\it arbitrary} flat $F$ will have at least the zero subspace $\{0\}$ in the
intersection $X_F \cap T$.
\
\item A {\it parabolic flat} $F$ is one for which $X_F \cap T$ is of maximum possible dimension, that is,
$\dim(X_F)$.  In this case, it must contain a cone $\sigma$ of this same dimension, say
indexed by the coset $wW_J$, whose linear span has pointwise stabilizer
$wW_Jw^{-1}$.  Hence $F = w\Phi^+_J$ is a $W$-conjugate\footnote{Strictly speaking, in order to
insure that $w\Phi^+_J \subseteq \Phi^+$, we should insist here that the coset representative
$w$ for $wW_J$ is chosen to be of minimum length.} of a {\it standard parabolic} flat $\Phi_J^+$,
where $\Phi_J^+$ is the subset of positive roots lying in the span of the simple roots $\{\alpha_s\}_{s \in J}$.

\item A {\it finite parabolic flat} $F$ is one for which $X_F \cap \interior(T)$ has maximum possible
dimension $\dim(X_F)$.  This turns out \cite[Cor. 3.8]{Walter}
to be equivalent to $F$ being the parabolic flat $w\Phi_J^+$
where the parabolic subgroup
$W_J$ is {\it finite}.  In other words,
$\interior(T)$ is the union of the cones $\sigma$ in $T$ whose associated parabolic subgroup is finite.
\end{enumerate}

It turns out that $W$ (or equivalently, $\Phi$) is finite if and only if the Tits cone $T$ coincides with
the whole space $V^*$ (and hence also coincides with the interior $\interior(T)$).
In this case, there is no distinction between the three kinds of flats:  all flats are
finite parabolic.  The reader interested solely in the case of finite Coxeter groups $W$ can therefore
safely ignore the remainder of this section.

\begin{example} \rm
The distinctions between the three kinds of flats are well-illustrated by
the case where $(W,S)$ is an irreducible {\it affine} Coxeter system, that is,
when the bilinear form $(\cdot,\cdot)$ is positive semidefinite, but degenerate.
In this case, the kernel of the bilinear form is a $1$-dimensional subspace $\ell$, and one
can faithfully realize the group $W$ as one generated by {\it Euclidean affine reflections}
in an affine hyperplane of $V^*$:  simply intersect the Tits cone with a (strictly) affine hyperplane
normal to $\ell$. See \cite[Ch. V, Sec. 4.9]{Bourbaki}.

As an example of a non-parabolic flat in this situation, pick any two roots $\alpha,\beta$
in $\Phi^+$ whose corresponding affine reflections $s_\alpha, s_\beta$ have parallel reflecting hyperplanes.
Then $s_\alpha, s_\beta$ generate an infinite subgroup $W' \subsetneq W$, whose
fixed subspace in $V^*$ corresponds to a flat $F=\cl(\{\alpha,\beta\})$
that is {\it not} parabolic.  To see this,
note that in the affine case, all proper parabolic subgroups
$W_J$ with $J \subsetneq S$ are finite, and hence all proper parabolic subgroups $wW_Jw^{-1}$
are also finite.  But $W'$ is infinite.

There is also a unique parabolic flat which is not finite parabolic in this situation,
namely the improper flat $F=\Phi^+$.  Its
corresponding intersection subspace $X_F=\{0\}$ lies in the Tits cone $T$, but not in its interior $\interior(T)$.
\end{example}

\begin{remark} \rm
The geometry of the Tits cone when $W$ is infinite, and in particular, its
interior $\interior(T)$, turn out to be important in the geometric group theory surrounding
the {\it generalized (Artin) braid group} $B(W,S)$ associated to $(W,S)$.
Inside the complex vector space $V^* \otimes {\mathbb C}$,
one has the open subset $V^* + i \cdot \interior(T)$, from which one can remove the intersection
with the complexified hyperplanes $\bigcup_{\alpha \in \Phi^+} H_\alpha$.  This hyperplane complement
carries a free action of $W$ with interesting topology:  it is conjectured
(and proven in many cases) that it is an Eilenberg-MacLane $K(PB(W,S),1)$-space
for the (pure) braid group $PB(W,S)$, and hence that its quotient by
$W$ is a $K(B(W,S),1)$ for the braid group $B(W,S)$; see Charney and Davis \cite{CharneyDavis}.
\end{remark}

In principle one might therefore consider (at least) three different versions  of the
poset of flats of the oriented matroid $M_\Phi$: the posets of arbitrary, parabolic, or
finite parabolic flats
$$
L^{\arbitrary}(M_\Phi) \supset  L^{\parabolic}(M_\Phi) \supset L^{\finiteparabolic}(M_\Phi).
$$
By the previous discussion, these posets of flats are isomorphic to
the following posets of subgroups.

\begin{proposition}
The map
$$
W' \mapsto \{\alpha \in \Phi^+: s_{\alpha} \in W'\}
$$
induces isomorphisms between the following
posets of subgroups and posets of flats,
all ordered by inclusion:
$$
\begin{aligned}
\{\text{ reflection subgroups }\} &\cong L^{\arbitrary}(M_\Phi) \\
\{\text{ parabolic subgroups }\} &\cong L^{\parabolic}(M_\Phi)\\
\{\text{ finite parabolic subgroups }\} &\cong L^{\finiteparabolic}(M_\Phi).
\end{aligned}
$$
\end{proposition}

When $W$ is finite, of course, all three notions of flats coincide and
all reflection subgroups are finite parabolic.  When $W$ is infinite we must make a choice
of which flats to consider.

\begin{remark} \rm
When discussing the De~Concini-Procesi wonderful compactifications of
Coxeter arrangements, Carr and Devadoss \cite{CD} chose to consider only
Coxeter systems $(W,S)$ which they call {\it simplicial}, namely
those in which every proper parabolic subgroup $W_J$ with $J \subsetneq S$
is finite, or equivalently, the simplicial decomposition of the Tits cone $T$ intersected
with the unit sphere in $V^*$ is a locally finite simplicial complex.
Such Coxeter systems include those which are finite, affine and hyperbolic.
When doing the wonderful compactification, they made
the natural choice of compactifying the complement of the
arrangement within the {\it interior} $\interior(T)$ of the Tits cone, after intersecting with the sphere.
This means that they only blew up along the {\it finite parabolic flats}, those
in  $L^{\finiteparabolic}(M_\Phi)$,
and avoided the problem of how to define blow-ups along
non-parabolic flats, where
the normal structure is not that of a finite hyperplane arrangement.
\end{remark}

We make a slightly different choice.  If one is not so concerned with the blow-ups themselves, but
rather with the truncations of the fundamental simplex $R$ which would tile the hypothetical blow-up,
then these truncated polytopes (the graph associahedra) are well-defined whether
or not the arrangement is locally finite.  In particular, we would like to consider graph-associahedra
associated to graphs $\Gamma$ for which $(W,S)$ is {\it not} simplicial, such as the complete graphs $\Gamma=K_n$
for $n \geq 4$.  For this reason we do not restrict ourselves
to the finite parabolic flats; instead we consider
all parabolic flats.

On the other hand, the relevant flats of the matroid which are relevant
for these truncations and blow-ups in the
wonderful compactification are those which intersect the Tits cone $T$ in full
dimension.  For this reason, when discussing the positive Bergman complex
$\B^+(M)$ in the next section,
we will consider only the poset of
parabolic flats $L^{\parabolic}(M_\Phi)$.

\begin{remark} \rm
It is not clear that we should expect good behavior from the Bergman complex
$\B(M_\Phi)$ or positive Bergman complex $\B^+(M_\Phi)$ defined with respect to
the lattice $L^{\arbitrary}(M_\Phi)$ of {\it arbitrary} flats, when 
$M_{\Phi}$ is infinite.
\end{remark}

\section{The positive Bergman complex of a Coxeter arrangement}\label{Cox}

In this section we prove that the positive Bergman complex of a
Coxeter arrangement of type $\Phi$ is dual to the
graph associahedron of type $\Phi$.  More precisely, both of these
objects are homeomorphic to spheres of the same dimension, and
their face posets are dual.

\begin{caution}
Throughout this section, whenever the Coxeter system $(W,S)$ with
root system $\Phi$ has $W$ (or equivalently, $\Phi$) infinite,
the word {\it flat} used in the connection with the
oriented matroid $M_\Phi$ will mean a {\it parabolic flat}, as discussed in the end of Section~\ref{Tits-cone}.
\end{caution}

\begin{proposition}
\label{Coxeter-flats}
Let $(W,S)$ be an arbitrary Coxeter system, with root system $\Phi$ and Coxeter diagram
$\Gamma$.
\begin{enumerate}
\item[(i)] Positive flats in the oriented matroid $M_\Phi$ correspond to subsets $J \subset S$,
that is, they are the standard parabolic flats $\Phi_J^+$.
\item[(ii)] Connected positive flats in the oriented matroid $M_\Phi$
correspond to subsets
$J \subset S$ such that the vertex-induced subgraph $\Gamma_J$ is connected, that is, to tubes
in $\Gamma$.
\item[(iii)] The simple roots $\Delta$ form a circuitous base for the matroid $M_\Phi$.
\item[(iv)] If $F \subset G$ are flats in $M_\Phi$ with $G$ connected, then
the matroid quotient $G/F$ is connected.
\end{enumerate}
\end{proposition}
\begin{proof}

\noindent
{\sf (i)}: The hyperplanes bounding the base region/tope $R$ are $\{H_{\alpha_s}: s \in S\}$,
so positive flats are those spanned by sets of the form $\{\alpha_s: s \in J\}$ for subsets
$J \subset S$.  As in the previous section, we denote such a positive flat by $\Phi_J^+$.

\noindent
{\sf (ii)}:  Let $J \subset S$ with subgraph $\Gamma_J$, and
consider its associated
positive flat $\Phi_J^+$.
The first assertion of Lemma~\ref{KacLemma}
shows that $\Phi_J^+$ will
not be connected if $\Gamma_J$ is disconnected.
To see this, represent the flat $\Phi_J^+$ by a matrix in which
the rows correspond to simple roots of $\Phi_J^+$, i.e. vertices of
$\Gamma_J$, and the columns
express each positive root in $\Phi_J^+$ as a combination of simple roots.
By permuting columns, one can obtain a matrix which is a block-direct
sum of two smaller matrices, and hence $\Phi_J^+$ will not be connected.

On the other hand,
if $\Gamma_J$ is connected,
then the second assertion of Lemma~\ref{KacLemma} shows that there
is a positive root $\alpha$ with $\supp \alpha = \Gamma_J$, and
consequently $\{\alpha_s: s \in J\} \cup \{\alpha\}$ gives a circuit in $M_\Phi$
spanning this flat, so it is connected.

\noindent
{\sf (iii)}: This follows from the argument in (ii);  given $J \subset S$ with
$\Gamma_J$ connected, the basic circuit $\circuit(\Delta,\alpha)$ where $\supp \alpha = \Gamma_J$
spans the connected flat corresponding to $J$.

\noindent
{\sf (iv):} Let the flats $F, G$ correspond (since they are assumed to be parabolic flats)
to the parabolic subgroups $uW_Ju^{-1}, vW_Kv^{-1}$.  Equivalently,
assume they are equal to $u \Phi^+_J, v \Phi^+_K$.  One can make the
following reductions:
\begin{enumerate}
\item[$\bullet$] Translating by $v^{-1}$, one can assume that  $v$ is the identity.
\item[$\bullet$] Since $(W_K, K)$ itself forms a Coxeter system with root system $\Phi_K$,
one can assume $M_\Phi=G$ and $K=S$.  In particular, $M_{\Phi}$ is
connected.
\item[$\bullet$] Replacing the Coxeter system $(W,S)$ by the system $(W,uSu^{-1})$,
one can assume that $u$ is the identity.
\end{enumerate}
In other words, $F$ is the positive flat corresponding to some
subgraph $\Gamma_J$ of $\Gamma$, and we must show $M_{\Phi}/F$ is
a connected matroid. This is a consequence of {\sf (iii)} and
Lemma \ref{lemma2}.
%
\end{proof}

We now
give our main result.

\def\thetheoremB{\ref{Theorem1}}
\begin{theoremB}
Let $(W,S)$ be an arbitrary Coxeter system, with root system $\Phi$,
Coxeter diagram $\Gamma$, and
associated oriented matroid $M_\Phi$.
Then the face
poset of the coarse subdivision of $\B^+(M_{\Phi})$ is dual to the
face poset of the graph associahedron $P(\Gamma)$.
\end{theoremB}

\begin{proof}

By Theorem \ref{Face}, we need to show that the face poset of
(the coarse subdivision of) $\B^+(M_{\Phi})$
is equal to the poset of tubings of $\Gamma$, ordered by
containment.  We begin by describing a map $\Psi$ from flags  of
positive flats to tubings of $\Gamma$.

By Proposition \ref{Coxeter-flats}, positive flats of $M_{\Phi}$
correspond to subsets $J \subset S$ or subgraphs $\Gamma_J$ of
the Coxeter graph $\Gamma$.  Furthermore, a positive flat is
connected if and only if $\Gamma_J$ is a tube, and hence an arbitrary
positive flat corresponds to a disjoint union of compatible tubes,
no two of which are nested.  Since an inclusion of flats corresponds
to an inclusion of the subsets $J$, a flag $\FF$ of positive flats corresponds
to a nested chain of such unions of non-nested compatible tubes, that is,
to a tubing $\Psi(\FF)$.
Furthermore, in this correspondence,
inclusion of flags corresponds to containment
of tubings.

We claim that the map from flags to tubings is surjective.  Given
some tubing of $\Gamma$, linearly order its tubes $J_1,\ldots,J_k$
by any linear extension of the inclusion partial ordering, and then the
flag $\FF$ of positive flats having $F_i$ spanned by
$\{\alpha_s: s \in J_1 \cup J_2 \cup \cdots \cup J_i \}$
will map to this tubing.

Lastly, we show that $\Psi$ is actually a well-defined injective map
when regarded as a map on cells of the coarse subdivision of
$\B^+(M_\Phi)$.  To do so, it is enough to show that
two flags $\FF, \FF'$ of positive flats
give the same tubing if and only if
$M_{\FF}$ and $M_{\FF'}$ coincide.
By Lemma~\ref{Coxeter-flats}(iv)
and Proposition~\ref{circuitous-prop}, we need to show that
$\Psi(\FF)$ and $\Psi(\FF')$ coincide if and only if
$T_\FF$ and $T_{\FF'}$ coincide.  But this is clear, because
by construction, the rooted forest $T_\FF$
ignores the ordering within the
flag, and only records the data of the tubes which appear,
that is, the tubing.
\end{proof}

\begin{corollary} The positive Bergman complex
of a Coxeter arrangement is simplicial, and is in fact,
a flag simplicial sphere.
\end{corollary}

Another corollary of our proof is a new realization for the
positive Bergman complex of a Coxeter arrangement: we can obtain
it from a simplex by a sequence of {\it stellar subdivisions}
(since stellar subdivisions are dual to the  
truncations defining the graph-associahedra; see \cite[Exercise 3.0]{Ziegler}).

\section{The Bergman complex of a Coxeter arrangement}\label{Nested}

In this section we will give a concrete description of the Bergman 
complex of a Coxeter arrangement, in terms of the nested set 
complex.  We will also address a question of 
Eugene Tevelev \cite{Tevelev} concerning the relationship of the
positive Bergman complex to the Bergman complex in this setup.

Nested set complexes are simplicial complexes at the combinatorial
heart of De~Concini and Procesi's subspace
arrangement models \cite{DP}, and of the resolution of
singularities in toric varieties \cite{FK}. We now recall the
definition of the minimal nested set complex of a meet-semilattice
$L$, which we will simply refer to as {\bf the} nested set complex
of $L$, and denote $\N (L)$.  For the sake of avoiding the technicalities
of infinite semilattices, matroids, and Coxeter groups, we will assume that
{\bf everything is finite} in this section.

Say an element $y$ of $L$ is \emph{irreducible} if the lower
interval $[\hat{0},y]$ cannot be decomposed as the product of
smaller intervals of the form $[\hat{0},x]$. The nested set
complex $\N (L)$ of $L$ is a simplicial complex whose vertices are
the irreducible elements of $L$. A set $X$ of irreducibles is
\emph{nested} if for any nonempty antichain $\{x_1, \ldots, x_k\}$ in $X$,
$x_1 \vee \cdots \vee x_k$ is {\it not} irreducible. These nested sets
are the simplices of $\N (L)$.

If $M$ is a matroid and $L_M$ is its lattice of flats, we will
also call $\N (L_M)$ the nested set complex of $M$, and denote it
$\N (M)$. It is easy to see that the irreducible elements of $L_M$ are the
connected flats of $M$.  

%

%
%

\def\thetheoremB{\ref{Theorem2}}
\begin{theoremB}
For any finite Coxeter system $(W,S)$ and associated finite root system $\Phi$,
the coarse subdivision of the Bergman complex $\B(M_{\Phi})$ of
the Coxeter arrangement of type $\Phi$ is equal to the nested set
complex $\N (M_{\Phi})$.  In particular, the Bergman complex $\B(M_{\Phi})$ is a
simplicial complex.
\end{theoremB}

We offer two proofs of this result.  The first is short, but not very self-contained
in that it invokes a result of Feichtner and Sturmfels.  They showed that, for any matroid $M$,
the simplicial complex $\N (M)$ has a geometric realization which is intermediate
in coarseness between the fine and coarse subdivisions of the Bergman complex
$\B(M)$.  Furthermore,
they gave the following criterion for when $\N(M)$ and the coarse subdivision of $\B(M)$ coincide.

\begin{theorem}\cite[Theorem 5.3]{FS}
The nested set complex $\N (M)$ and the Bergman 
complex $\B(M)$ together with its coarse subdivision
coincide if and only if the matroid $G/F$ is connected for every pair of flats $F \subset G$ 
in which $G$ is connected.
\end{theorem}

\noindent
Theorem~\ref{Theorem2} then follows immediately from this result together
with our Proposition~\ref{Coxeter-flats}(iv).

\vskip.2in

  On the other hand, one might suspect that in the case of a Coxeter arrangement,
Theorem~\ref{Theorem1} describing the positive Bergman complex $\B^+(M)$,
and Theorem~\ref{Theorem2} about the entire Bergman complex $\B(M)$,
should be related by Theorem~\ref{posBergmanscoverBergman} 
showing that $\B(M)$ is covered by several copies of $\B^+(M)$.  The goal of the
remainder of this section is to develop this connection, in the context of arbitrary
oriented matroids, partly in order to answer the question of Tevelev mentioned earlier,
and partly for its own interest.

  We begin in the setting of an (unoriented) matroid $M$, giving the relationship between
nested sets and the labelled forest $T_\FF$ associated to a flag of flats of 
$M$ 
in Proposition \ref{forest-determines-matroid}.  This next proposition can be gleaned
implicitly from the material in \cite[Sections 3 and 4]{FS}, but we state it explicitly
here, and include our own proof, for the sake of self-containment.

\begin{proposition}
\label{nested-sets-are-forests}
Let $M$ be an (unoriented) matroid.  
Given a flag $\F$ of flats in $M$, consider the set of
connected flats $\{G_i\}$ which label the vertices of the forest $T_\FF$.
Then the collection of all such sets $\{G_i\}$ of flats, as 
$\F$ ranges over all flags of flats in $M$, are precisely the nested sets of $M$.
\end{proposition}

\begin{proof}
For any $\F$, the labels of $T_{\F}$ are connected flats of $M$ by definition. 
Now let us show that they form a nested set. We need to show that for any antichain 
$\{G_1, \ldots, G_k\}$ of flats among the vertex labels $T_{\F}$, 
their $L_M$-join is not connected, so assume that it is. 
Let $F_i$ be the smallest flat of $\F$ containing $H=G_1 \vee \cdots \vee G_k$. 
Since $H$ is connected, it is a subset of a connected component $G$ of $F_i$. 
But $G$ lies
above $G_1, \ldots, G_k$ in the tree; this means that the $G_i$s are 
connected components of flats $F_j$ with $j<i$, and are therefore contained in 
$F_{i-1}$. But then $H$, being their $L_M$-join, must also be contained in 
$F_{i-1}$, contradicting the minimality of $i$.

Now we show that every nested set $N$ of $M$, when its elements are
ordered by inclusion, 
can be obtained as the vertex labels of the forest $T_\FF$ for some flag $\FF$ of 
flats of $M$. The connected flats in $N$ can be labelled $G_1, G_2, \ldots, G_k$ 
in such a way that $i<j$ implies $G_i \not\supseteq G_j$. 
Let the flag $\F$ consist of the flats $F_i = G_1 \vee  \cdots \vee G_i$ for $1 \leq i \leq k$. 

First we show that $F_i$ is just the union of the maximal $G_j$s with $j \leq i$. 
Suppose there was an element $e$ in $F_i$ which is not in one of these $G_j$s. 
Then there must be a circuit containing $e$ and elements of, say, $G_a, \ldots, G_z$. 
Then $e, G_a, \ldots, G_z$ are contained in the same connected component of $M$. 
Any other $f$ in $G_a \vee \cdots \vee G_z$ is in this same component, 
so $G_a \vee \cdots \vee G_z$ is connected. This contradicts the assumption that $N$ is nested.

By the same reasoning, we cannot have a circuit in $F_i$ consisting of elements of 
more than one of the $G_j$s. Therefore the maximal $G_j$s with $j \leq i$ are 
actually the connected components of $F_i$. In particular, $G_i$ is one of them. 
This shows that the flag $\F$ gives rise to the nested set $N$, as we wished to prove.
\end{proof}

We next wish to understand, in the setting of an (acyclically) {\it oriented} matroid $M$, how 
nested sets interact with the notion of positive flats.  First, we need a small technical lemma.

\begin{lemma}
\label{components-of-positive-flats}
If a flat of an oriented matroid $M$ is positive, so are its connected components.
\end{lemma}

\begin{proof}
 Let $G$ be a connected component of a positive flat $F$ and assume for the sake of contradiction
that $G$ is not positive. By \cite[Proposition 9.1.2]{RedBook}, we can find a signed circuit 
$X$ of $M$ such that $X^+ \subseteq G$ and $X^- \not\subseteq G$. 
We then have that $X^+ \subseteq F$, which implies that $X^- \subseteq F$ 
since $F$ is positive. Therefore $X$ is a circuit in $F$, containing elements of 
more than one of its connected components. This is a contradiction.
\end{proof}

\begin{proposition}
\label{oriented-nested-sets-are-forests}
Let $M$ be an acyclically oriented matroid whose positive tope is simplicial. 
As $\F$ ranges over all flags of positive flats in $M$, the sets of flats labelling vertices of
the forests $T_{\F}$ are precisely those nested sets of $M$ 
which consist of positive (connected) flats.
\end{proposition}

\begin{proof}
By Proposition \ref{nested-sets-are-forests}, if $\F$ is a flag of positive flats, 
then $T_{\F}$ is a nested set of $M$. The labels of $T_{\F}$ are positive by 
Lemma \ref{components-of-positive-flats}. 

Now start with a nested set $N$ of $M$ consisting of positive connected flats, 
labelled $G_1, G_2, \ldots, G_k$ in such a way that $i<j$ implies $G_i \not\supseteq G_j$. 
As in the proof of Proposition \ref{nested-sets-are-forests}, the 
flag $\F$ consisting of the flats $F_i = G_1 \vee  \cdots \vee G_i$ 
for $1 \leq i \leq k$ satisfies $T_{\F} = N$. Finally, each $F_i$ is a disjoint union of $G_j$s, 
which are positive. Since the positive tope is simplicial, the $F_i$s are also positive.
\end{proof}

\begin{remark} \rm \ \\
Proposition~\ref{oriented-nested-sets-are-forests} is closely related to 
Feichtner and Sturmfels'
notion \cite[Section 4]{FS} of the {\it localization} of the nested set complex $\N (M)$ to a basis $B$ of the matroid
$M$, if one chooses $B$ to be the elements of the ground set which bound the simplicial positive tope 
of $M$.
\end{remark}

\begin{observation} \rm \ \\
Proposition \ref{oriented-nested-sets-are-forests} can fail if the positive tope is not simplicial. For example, consider the oriented matroid $M$ of affine dependencies of the vertices of a square which are cyclically labelled
$\{1, 2, 3, 4\}$. Now, $\{1,3\}$ is a nested set of $L_M$ consisting of positive flats. However, it does not arise as the forest of a flag of positive flats.
\end{observation}

\begin{observation} \rm \ \\
The nested sets of $M$ which consist of positive (connected) flats are not the same as the nested sets of the lattice of positive flats $\F_{lv}(M)$. In the previous example, $\{1,3\}$ is nested in $L_M$ but not in $\F_{lv}(M)$. Even in Example \ref{braid-example}, the (simplicial) oriented matroid of the braid arrangement $A_3$, $\{1,4\}$ is not nested in $L_M$ but it is nested in $\F_{lv}(M)$.
\end{observation}


We now give the second proof of Theorem \ref{Theorem2}.
In view of Proposition \ref{nested-sets-are-forests}, 
Proposition \ref{forest-determines-matroid} tells us that, for any matroid $M$, 
the nested set complex $\N(M)$ is a refinement of the coarse subdivision of
the Bergman complex $\B(M)$ and a 
coarsening of the order complex $\Delta(\overline{L}_M)$ (i.e.\ the fine 
subdivision of $\B(M)$). 
Therefore, it suffices to show that, in the case of a (finite) Coxeter arrangement, 
every cell in the nested set complex is equal to the cell of the Bergman complex which contains it.

So consider an arbitrary cell $C_N$ corresponding to a nested set $N$ in the nested set complex 
$\N(M_{\Phi})$, and the cell $D$ of the Bergman complex $\B(M_{\Phi})$ containing it. 
By Theorem 2.8, we can find a tope $T$ of $M$ such that the positive Bergman complex $\B^+(M_{\Phi})$ 
corresponding to $T$ contains the cell $D$, and therefore the cell $C_N$. 
Since $M_{\Phi}$ is simplicial, Proposition \ref{oriented-nested-sets-are-forests} implies 
that the flats in $N$ are positive with respect to $T$. 

Proposition \ref{Coxeter-flats} tells us that connected positive flats correspond to tubes in $\Gamma$. 
It is easy to see that a set of connected positive flats is nested if and only if the corresponding set of 
tubes is a tubing of $\Gamma$. Therefore $C_N$ is precisely the cell of $\B^+(M_{\Phi})$ labelled by the 
tubing corresponding to $N$, by (the proof of) Theorem \ref{Theorem1}. 
It follows that $C_N=D$, as we wished to show.
$\qed$

\vskip.2in

Recently Tevelev \cite{Tevelev} asked whether every (coarse) cell in
the Bergman complex $\B(M_{\Phi})$ of a Coxeter arrangement of type
$\Phi$ is Coxeter-group equivalent to a 
cell in $\B^+(M_{\Phi})$, i.e. a cell obtained from a tubing
of the corresponding Dynkin diagram.  Using Theorem \ref{posBergmanscoverBergman}, we 
can give an affirmative answer to this question.  We begin with the following observation
about the group $\Aut(M)$ of all automorphisms 
$\phi: E \rightarrow E$ for a finite matroid $M$ on ground set $E=[n]$.

\begin{proposition} \label{action}
The group $\Aut(M)$, acting on $\R^n$ by permuting coordinates, preserves the Bergman fan $\B(M)$, and
acts by cellular automorphisms on its coarse subdivision.
\end{proposition}

\begin{proof}
Recall that a weight vector $\omega$ in $\R^n$ has a uniquely associated flag $\FF(\omega)$ 
of subsets as in \eqref{flag-of-flats}, and any permutation $\phi: E \rightarrow E$ respects this association:
$\FF(\phi(\omega))=\phi(\FF(\omega))$.  Since $\omega$ lies in $\B(M)$ if and only if
this flag of subsets is a flag of flats, a matroid automorphism $\phi$ will preserve this property,
and hence preserves $\B(M)$.

Recall also from \eqref{matroid-from-flag} that the matroid $M_{\omega}$ induced by $\omega$
is exactly $\bigoplus_{i=1}^{k+1} F_i / F_{i-1}
$.  Since
two weight vectors $\omega, \nu$ lie in the same coarse cell if
and only if $M_{\omega} = M_{\nu}$, the second
assertion of the proposition follows.
\end{proof}

\begin{proposition}
Let $M_{\Phi}$ be the oriented matroid of a (finite) Coxeter arrangement
of type $\Phi$ with Coxeter group $W$.  Then any
coarse cell in $\B(M_{\Phi})$ is $W$-equivalent to 
a coarse cell in $\B^+(M_{\Phi})$.
\end{proposition}

\begin{proof}
Since $W$ acts by matroid automorphisms on
the ground set $E$ of $M_\Phi$, Proposition~\ref{action}
implies that $W$ permutes the coarse cells of $\B(M)$.

Now let $C$ be any coarse cell of $\B(M)$, and choose a fine 
cell $c \subset C$ defined by a flag of flats $\FF$.  By
Theorem \ref{posBergmanscoverBergman} and the fact that $W$ acts
transitively on the regions of the arrangement, 
there exists some $w \in W$ such that $w(c)$ is a fine cell
defined by a flag of positive flats, i.e. $w(c)$ lies inside a 
coarse cell $D$ of $\B^+(M)$.  But now it follows that 
$w(C) = D$: $w(C)$ must be a cell of $\B(M)$, and it contains
$w(c)$, which lies in $D$. 
\end{proof}


\medskip

\textsc{Acknowledgments} We are very grateful to Satyan Devadoss for
allowing us to reproduce two of his figures from \cite{CD}.  These are
our Figures \ref{D4} and \ref{A3}.  Additionally, we thank Bernd
Sturmfels, Eva Feichtner, and Eugene Tevelev for useful discussions.

\raggedright


\addcontentsline{toc}{section}{References}
\begin{thebibliography}{2}


\bibitem{Berg} F. Ardila and C. Klivans,
\emph{The Bergman complex of a matroid and phylogenetic trees}, 
to appear in the Journal of Combinatorial Theory, Series B.
\textsf{arXiv:math.CO/0311370}

\bibitem{PosBerg} F. Ardila, C. Klivans, and L. Williams,
\emph{The positive Bergman complex of an oriented matroid}, to
appear in the European Journal of Combinatorics.


\bibitem{BjornerBrenti} A. Bj\"orner and F. Brenti,
\emph{The combinatorics of Coxeter groups},
Graduate Texts in Mathematics, {\bf 231}.
Springer, New York, 2005.

\bibitem{RedBook}
A. Bj\"{o}rner, M. Las Vergnas, B. Sturmfels, N. White, G. 
Ziegler, \emph{Oriented matroids,}  Encyclopedia of Mathematics
and its Applications, vol. 46,  Cambridge University Press,
Cambridge, 1993.

\bibitem{BW}
A. Bj\"{o}rner and M. Wachs, \emph{Geometrically constructed bases
for homology of partition lattices of types $A$, $B$ and $D$,}
Electronic Journal of Combinatorics {\bf 11 (2)} (2004).

\bibitem{Bourbaki}
N. Bourbaki,
\emph{Lie groups and Lie algebras} Chapters 4--6,
Elements of Mathematics (Berlin),
Springer-Verlag, Berlin, 2002.

\bibitem{Brown}
K.S. Brown,
\emph{Buildings},
Springer Monographs in Mathematics, Springer-Verlag, New York, 1998.

\bibitem{CD}
M. Carr, S. Devadoss,  \emph{Coxeter complexes and graph associahedra},
to appear in Topology and its Applications.

\bibitem{CharneyDavis}
R. Charney and M.W. Davis,
\emph{Finite $K(\pi, 1)$s for Artin groups},
Prospects in topology (Princeton, NJ, 1994),  110--124,
Ann. of Math. Stud. {\bf 138},
Princeton Univ. Press, Princeton, NJ, 1995.

\bibitem{Davis}
M. Davis, T. Januszkiewicz, and R. Scott, \emph{Fundamental groups of minimal
blow-ups}, Adv. Math. {\bf 177} (2003) 115-179.

\bibitem{DP}
C. De~Concini, C. Procesi,  \emph{Wonderful models of subspace
arrangements}, Selecta Math. (N.S.) 1, no. 3 (1995), 459-494.

\bibitem{FK} E. Feichtner, D. Kozlov,
\emph{Incidence combinatorics of resolutions},
Selecta Math. (N.S.) 10 (2004), 37-60.

\bibitem{FM} E. Feichtner, I. M\"{u}ller, \emph{On the topology of
nested set complexes},  Proc. Amer. Math. Soc. 133 (2005), 999-1006.

\bibitem{FS} E. Feichtner, B. Sturmfels,
\emph{Matroid polytopes, nested sets and Bergman fans}, to appear in
\emph{Port. Math.}.

\bibitem{GZ} C. Greene, T. Zaslavsky, \emph{On the interpretation of Whitney numbers
through arrangements of hyperplanes, zonotopes, non-Radon
partitions, and orientations of graphs}, Trans. Amer. Math. Soc.
{\bf 280} (1983), 97-126.

\bibitem{Humphreys} J. Humphreys,
\emph{Reflection groups and Coxeter groups.} Cambridge University Press,
1990.

\bibitem{Kac}
V. Kac, \emph{Infinite dimensional Lie algebras,}  Cambridge University
Press, 1990.

\bibitem{LasVergnas}
M. Las Vergnas, \emph{Acyclic and totally cyclic orientations of combinatorial geometries}, Discrete Math. {\bf 20} (1977), 51-61.

\bibitem{Postnikov} A. Postnikov, \emph{Permutahedra, associahedra,
and beyond}, preprint, 2005.
\textsf{arXiv:math.CO/0507163}

\bibitem{Speyer}
D. Speyer and B. Sturmfels, \emph{The tropical Grassmannian,}
Adv. Geom., {\bf 4} (2004), 389-411.

\bibitem{SpeyerWilliams}
D. Speyer and L. Williams, \emph{The tropical totally positive
Grassmannian}, Journal of Algebraic Combinatorics 22 (2005), 189-210.

\bibitem{EC1}
R. P. Stanley, \emph{Enumerative combinatorics, vol.1,} Cambridge
University Press, New York, 1986.

\bibitem{Sturmfels}
B. Sturmfels, \emph{Solving systems of polynomial equations,} CBMS
Regional Conference Series in Mathematics, vol. 97,  American
Mathematical Society, Providence, RI, 2002.
\bibitem{Tevelev}
E. Tevelev, personal communication

\bibitem{Walter}
J.H. Walter,
\emph{Bases of chambers of linear Coxeter groups},
Groups and combinatorics---in memory of Michio Suzuki,
423--436, {\emph Adv. Stud. Pure Math.} {\bf 32},
Math. Soc. Japan, Tokyo, 2001.


\bibitem{Zaslavksy}
T. Zaslavsky,
\emph{The M\"obius function and the characteristic polynomial},
Chapter 7 in Theory of matroids (edited by N. White)
Encyclopedia of Mathematics and its Applications, {\bf 26}.
Cambridge University Press, Cambridge, 1986.

\bibitem{Ziegler}
G. Ziegler, \emph{Lectures on polytopes}, Graduate Texts in Mathematics,
vol. 152, Springer-Verlag, New York, 1995.

\end{thebibliography}
\end{document}